\documentclass[leqno,12pt]{article}

\usepackage{amssymb}
\usepackage{mathrsfs} 
\usepackage{amsmath}
\usepackage{amsthm}
\usepackage{amsfonts}
\usepackage{verbatim} 
\usepackage{graphicx}

\def\thesection{\arabic{section}}
\renewcommand{\theequation}{\thesection.\arabic{equation}}

\newtheorem{theorem}{Theorem}[section]
\newtheorem{lemma}[theorem]{Lemma}

\newtheorem{corollary}[theorem]{Corollary}

\newtheorem{definition}[theorem]{Definition}

\theoremstyle{definition}   
\newtheorem{remark}[theorem]{Remark}

\evensidemargin\oddsidemargin

\newcommand{\eqnsection}{
\renewcommand{\theequation}{\thesection.\arabic{equation}}
    \makeatletter
    \csname  @addtoreset\endcsname{equation}{section}
    \makeatother}
\eqnsection

\def\P{{\bf P}}
\def\E{{\bf E}}

\def\z{{\mathbb Z}}

\def\ee{\mathrm{e}}
\def\d{\, \mathrm{d}}

\begin{document}

\baselineskip=18pt
\setcounter{page}{1}


\vglue50pt

\centerline{\large\bf The stable Derrida--Retaux system at criticality}

\bigskip
\bigskip

\centerline{\footnotesize Xinxing Chen\footnote{\scriptsize School of Mathematical Sciences, Shanghai Jiaotong University, 200240 Shanghai, China, {\tt chenxinx@sjtu.edu.cn} $\;$Partially supported by NSFC grant Nos.  11771286 and 11531001.}
and
Zhan Shi\footnote{\scriptsize LPSM, Sorbonne Universit\'e Paris VI, 4 place Jussieu, F-75252 Paris Cedex 05, France, {\tt zhan.shi@upmc.fr} $\;$Partially supported by ANR MALIN.}
} 

\bigskip
\bigskip

\centerline{\it Dedicated to the memory of Vladas Sidoravicius}

\bigskip
\bigskip

{\leftskip=1.5truecm \rightskip=1.5truecm \baselineskip=15pt \small

\noindent{\slshape\bfseries Summary.} The Derrida--Retaux recursive system was investigated by Derrida and Retaux~\cite{derrida-retaux} as a hierarchical renormalization model in statistical physics. A prediction of \cite{derrida-retaux} on the free energy has recently been rigorously proved (\cite{bmvxyz_conjecture_DR}), confirming the Berezinskii--Kosterlitz--Thouless-type phase transition in the system. Interestingly, it has been established in \cite{bmvxyz_conjecture_DR} that the prediction is valid only under a certain integrability assumption on the initial distribution, and a new type of universality result has been shown when this integrability assumption is not satisfied. We present a unified approach for systems satisfying a certain domination condition, and give an upper bound for derivatives of all orders of the moment generating function. When the integrability assumption is not satisfied, our result allows to identify the large-time order of magnitude of the product of the moment generating functions at criticality, confirming and completing a previous result in \cite{collet-eckmann-glaser-martin}.

\bigskip


\bigskip


} 

\bigskip
\bigskip

\section{Introduction}
\label{s:intro}

Fix an integer $m\ge 2$. Let $X_0$ be a random variable taking values in $\z_+ := \{ 0, \, 1, \, 2, \ldots\}$. To avoid trivial discussion, it is assumed, throughout the paper, that $\P ( X_0 \ge 2) > 0$. Let us consider the Derrida--Retaux recursive system $(X_n, \, n\ge 0)$ defined as follows: for all $n\ge 0$,
\begin{equation}
    X_{n+1}
    =
    (X_{n,1} + \cdots + X_{n,m} -1)^+ ,
    \label{iteration}
\end{equation}

\noindent where $X_{n,i}$, $i\ge 1$, are independent copies of $X_n$. This was investigated by Derrida and Retaux~\cite{derrida-retaux} as a toy model to study depinning in presence of impurities \cite{luck,derrida-hakim-vannimenus,tang-chate,giacomin,giacomin-toninelli,dglt,monthus}. We refer to \cite{bz_dobrushin} for an overview on rigorous results and predictions about the Derrida--Retaux system.

Assuming $\E(X_0)<\infty$, it is immediate from \eqref{iteration} that $\E(X_{n+1}) \le m\, \E(X_n)$, so the free energy
$$
F_\infty := \lim_{n\to \infty} \downarrow \frac{\E(X_n)}{m^n} \in [0, \, \infty),
$$

\noindent is well-defined. A remarkable result by Collet et al.~\cite{collet-eckmann-glaser-martin} tells us that assuming $\E(X_0 \, m^{X_0})<\infty$ (which we take for granted throughout the paper) and writing $\eta := (m-1) \E(X_0 \, m^{X_0}) - \E(m^{X_0})$, then $F_\infty >0$ if $\eta >0$, and $F_\infty =0$ if $\eta \le 0$.

As such, it is natural to say that the system $(X_n, \, n\ge 0)$ is supercritical if $\eta>0$, is critical if $\eta=0$, and is subcritical if $\eta<0$.

It has been conjectured by Derrida and Retaux~\cite{derrida-retaux} that if $\eta >0$, then we would have
\begin{equation}
    F_\infty = \exp \Big( - \frac{C+o(1)}{\eta^{1/2}} \Big), \qquad \eta \to 0+ \, ,
    \label{conj_DR}
\end{equation}

\noindent for some constant $C\in (0, \, \infty)$ possibly depending on the law of $X_0$. A (somehow weak) result has been proved in \cite{bmvxyz_conjecture_DR}: assuming $\E(X_0^3 m^{X_0})<\infty$,
$$
F_\infty = \exp \Big( - \frac{1}{\eta^{1/2+o(1)}} \Big), \qquad \eta \to 0+ \, .
$$

\noindent This confirms that the Derrida--Retaux system has a Berezinskii--Kosterlitz--Thouless-type phase transition of infinite order. The integrability assumption $\E(X_0^3 m^{X_0})<\infty$ might look exotic, but it is optimal. [We believe that there should be a change-of-measures argument, and that the assumption is equivalent to saying that $X_0$ has a finite second moment under a new probability measure; however, we have not succeeded in making this idea into a rigorous argument.] In fact, it has also been proved in \cite{bmvxyz_conjecture_DR} that if $\P(X_0 = k) \sim c \, m^{-k} k^{-\alpha}$, $k\to \infty$, for some $2<\alpha<4$ and $c>0$,\footnote{Notation: $a_k \sim b_k$, $k\to \infty$, means $\lim_{k\to \infty} \frac{a_k}{b_k} =1$.} then
\begin{equation}
    F_\infty = \exp \Big( - \frac{1}{\eta^{\nu +o(1)}} \Big), \qquad \eta \to 0+ \, ,
    \label{energie_cas_stable}
\end{equation}

\noindent where $\nu = \nu(\alpha) := \frac{1}{\alpha-2}$. In other words, \eqref{conj_DR} predicts only a small part of universalities, under the assumption $\E(X_0^3 m^{X_0})<\infty$, while other universality phenomena are described by \eqref{energie_cas_stable}. We expect many other universality results in the latter setting (for example, corresponding to those in \cite{4authors} for an analogous continuous-time model); unfortunately, they are currently only on a heuristic level.

It is well-known that sum of i.i.d.\ random variables, after an appropriate normalization, converges to a Gaussian limiting law under the condition of finiteness of second moment, and to a stable limiting law under a weaker integrability condition. We say that the Derrida--Retaux system has a ``finite variance" if $\E(X_0^3 m^{X_0})<\infty$, and that it is a stable system if integrability condition holds for lower orders. In this paper, we are interested in the stable system when it is critical, i.e., when $(m-1) \E(X_0 \, m^{X_0}) = \E(m^{X_0})$. [We are going to see in Section \ref{s:domination}, quite easily, that this implies $(m-1) \E(X_n \, m^{X_n}) = \E(m^{X_n})$ for all $n\ge 0$.] We write $(Y_n, \, n\ge 0)$ instead of $(X_n, \, n\ge 0)$ in order to insist on criticality. {F}rom now on, we assume $(Y_n, \, n\ge 0)$ to be a Derrida--Retaux system satisfying $(m-1) \E(Y_0 \, m^{Y_0}) = \E(m^{Y_0}) <\infty$, such that
\begin{equation}
    \P(Y_0 = k) \sim c_0 \, m^{-k} k^{-\alpha},
    \qquad
    k\to \infty ,
    \label{alpha}
\end{equation}

\noindent for some $2<\alpha<4$ and $c_0>0$. We intend to prove the following result.

\medskip

\begin{theorem}
\label{t:main}

 Let $(Y_n, \, n\ge 0)$ be such that $(m-1) \E(Y_0 \, m^{Y_0}) = \E(m^{Y_0}) <\infty$. Under assumption \eqref{alpha}, there exist constants $c_2\ge c_1 >0$ such that for all $n\ge 1$,
 $$
 c_1 \, n^{\alpha-2} \le \prod_{i=0}^{n-1} [\E(m^{Y_i})]^{m-1} \le c_2 \, n^{\alpha-2} .
 $$

\end{theorem}

\medskip

When the system is of ``finite-variance" (i.e., $\E(Y_0^3 m^{Y_0})<\infty$), the analogue of Theorem \ref{t:main} was known (\cite{collet-eckmann-glaser-martin}, \cite{bmxyz_questions}), and has played an important role in the study of the asymptotics of $\P(Y_n >0)$ and $\E(Y_n)$ in \cite{bmvxyz_sustainability}. It would be tempting to believe that Theorem \ref{t:main} could play an equally important role in the study of the same problems for the stable system.

Just like the usual random walk has a nice continuous-time analogue which is Brownian motion, the Derrida--Retaux system has analogues in continuous time (Derrida and Retaux~\cite{derrida-retaux}, Hu, Mallein and Pain~\cite{HMP}), defined via appropriate integro-differential equations. For the continuous-time analogue of the stable Derrida--Retaux system, see~\cite{4authors}. These continuous-time models have been studied in depth in \cite{derrida-retaux}, \cite{HMP} and \cite{4authors}, while most of the corresponding problems remain open for the original Derrida--Retaux system.

With the exception of the case $\alpha=3$, Theorem \ref{t:main} was already stated in Collet et al.~\cite{collet-eckmann-glaser-martin}: its proof in case $2<\alpha<3$ was indicated, whereas the proof in case $3<\alpha<4$ was only summarized in a ``very succinct account". By means of the notion of dominability (see the forthcoming Definition \ref{d:dominable}), we give a unified approach to the system in both situations, i.e., either it is stable (no need for discussions separately on the cases $2<\alpha<3$ and $3<\alpha<4$) or is of ``finite variance". Concretely, in both situations, we use a truncating argument by considering a bounded random variable defined by
$$
Z_0 = Z_0(M) := Y_0 \, {\bf 1}_{\{ Y_0 \le a(M)\} } ,
$$

\noindent where $a(M) \in [1, \, \infty]$ can be possibly infinite (in which case there is no need for truncation), whose value depends on an integer parameter $M\ge 1$. Consider the Derrida--Retaux system $(Z_n, \, n\ge 0)$ whose initial distribution is given by $Z_0$.\footnote{Strictly speaking, it is a {\it sequence} of Derrida--Retaux systems, indexed by $M$.} We prove, in Theorem \ref{t:dominable}, that in both situations, it is possible to choose a convenient value of $a(M)$ such that the new system $(Z_n, \, n\ge 0)$ is dominable (in the sense of Definition \ref{d:dominable}), while it is possible to connect the moment generating functions of $Y_n$ and $Z_n$. In Theorem \ref{t:H_nkupper}, we give an upper bound for the moment generating function of any dominable system $(Z_n, \, n\ge 0)$. As such, a combined application of Theorems \ref{t:dominable} and \ref{t:H_nkupper} will yield information for the moment generating function of the original Derrida--Retaux system, in both situations. In the stable case, it will yield Theorem \ref{t:main}, whereas in the case of ``finite variance", under a stronger integrability assumption on the law of $Y_0$, it will give the following result:

\medskip

\begin{theorem}
 \label{t:moment_ub}

 Let $(Y_n, \, n\ge 0)$ be such that $(m-1) \E(Y_0 \, m^{Y_0}) = \E(m^{Y_0})$. If $\E(s^{Y_0})<\infty$ for some $s>m$, then there exists a constant $c_3 >0$ such that for all integers $n\ge 1$ and $k \ge 1$,
 \begin{equation}
     \frac{\mathrm{d}^k}{\mathrm{d} u^k}\, \E(u^{Y_n}) \Big|_{u=m}
     \le
     k!\, \ee^{c_3 k} \, n^{k-1} \, .
     \label{eq_t:moment_ub}
 \end{equation}

\end{theorem}

\medskip

In the ``finite-variance" case $\E(Y_0^3 m^{Y_0})<\infty$, \eqref{eq_t:moment_ub} for $k\in \{ 1, \, 2, \, 3\}$ was known: the case $k=1$ is simple because by criticality,  $\E(Y_n m^{Y_n-1}) = \frac{1}{m(m-1)} \, \E(m^{Y_n})$ which is bounded in $n$ (\cite{collet-eckmann-glaser-martin}, \cite{bmxyz_questions}), the case $k=3$ was proved in \cite{bmxyz_questions}, and the case $k=2$, stated in \cite{bmvxyz_conjecture_DR}, follows immediately from the cases $k=1$ and $k=3$ by means of the Cauchy--Schwarz inequality. More generally, if $\E(Y_0^\ell m^{Y_0})<\infty$ for some integer $\ell \ge 1$, then for all $k\in [1, \, \ell] \cap \z$, it is quite easy to prove (\cite{bmvxyz_sustainability}) by induction in $k$, using the recursion \eqref{iteration}, that there exists a constant $c>0$ such that for all integer $n\ge 1$,
$$
\frac{\mathrm{d}^k}{\mathrm{d} u^k}\, \E(u^{Y_n}) \Big|_{u=m}
\le
c \, n^{k-1} \, .
$$

\noindent Theorem \ref{t:moment_ub} gives information about the dependence in $k$ of the constant $c$, under the integrability assumption $\E(s^{Y_0})<\infty$ for some $s>m$.

The rest of the paper is organized as follows. In Section \ref{s:domination}, we introduce the notion of dominable systems. Theorem \ref{t:H_nkupper}, which gives an upper bound for the moment generating function of dominable systems, is the main technical result of the paper. The brief Section \ref{s:moment_s>m} is devoted to the proof of Theorem \ref{t:moment_ub}, obtained as a simple consequence of Theorem \ref{t:H_nkupper}. In Section \ref{s:truncating}, for both ``finite-variance" and stable systems, we construct a dominable system $(Z_n, \, n\ge 0)$ such that $Z_0$ is obtained from an appropriate truncation of $Y_0$. Finally, Theorem \ref{t:main} is proved in Section \ref{s:pf_thm}, also as a consequence of Theorem \ref{t:H_nkupper}.

\section{Dominable systems}
\label{s:domination}

We introduce the notion of dominable systems and prove a general upper bound for the moment generating function of such systems (Theorem \ref{t:H_nkupper}). As before, we talk about the Derrida--Retaux system $(Z_n, \, n\ge 0)$, while it is, in fact, a sequence of Derrida--Retaux systems $(Z_n(M), \, n\ge 0)$ indexed by the integer-valued parameter $M$.

\medskip

\begin{definition}
\label{d:dominable}

 Let $\gamma>0$. The system $(Z_n, \, n\ge 0)$ is said to be $\gamma$-dominable if for all sufficiently large integer $M$, say $M\ge M_0$, $Z_0=Z_0(M)$ is bounded, and there exists a constant $\vartheta(M) \ge 1$ such that $M \mapsto \vartheta(M)$ is non-decreasing in $M\ge M_0$, and that$\,$\footnote{In \eqref{assump_H0k} and \eqref{assump_Qn}, $k\ge 3$ and $n\ge 1$ are integers. Notation: $a\vee b := \max\{ a, \, b\}$.}
\begin{eqnarray}
    \E(Z_0^k m^{Z_0})
 &\le& M^{k-3} (\vartheta(M) + k!) ,
    \qquad
    k\ge 3,
    \label{assump_H0k}
    \\
    \vartheta (n\vee M) \prod_{i=0}^{n-1} [\E(m^{Z_i})]^{m-1}
 &\le& \gamma \, (n\vee M)^2 ,
    \qquad
    n\ge 1.
    \label{assump_Qn}
\end{eqnarray}

\end{definition}

\medskip

\begin{remark}
\label{r:dominable=>subcritical}

Let $(Z_n, \, n\ge 0)$ be $\gamma$-dominable. By \eqref{assump_Qn} and the trivial inequality $\E(m^{Z_i}) \ge 1$ ($\forall i\ge 0$), we have $\E(Z_n) \le \E(m^{Z_n}) \le [\gamma \, (n\vee M)^2]^{1/(m-1)}$, so the free energy $F_\infty := \lim_{n\to \infty} \frac{\E(Z_n)}{m^n}$ vanishes; by the criterion of Collet et al.~\cite{collet-eckmann-glaser-martin} recalled in the introduction, the system $(Z_n, \, n\ge 0)$ is subcritical or critical for all $M\ge M_0$: we have $(m-1)\E(Z_0 \, m^{Z_0}) \le \E(m^{Z_0})$.\qed

\end{remark}

\medskip

\begin{theorem}
 \label{t:H_nkupper}

 Let $\gamma>0$. Let $(Z_n, \, n\ge 0)$ be a $\gamma$-dominable system, and let $H_n(u) := \E(u^{Z_n})$. There exists a constant $c_4\ge 1$, depending only on $(m, \, \gamma)$, such that for all integers $k\ge 3$, $n\ge 1$ and $M\ge M_0$,
 $$
 H_n^{(k)}(m)\le k!\, c_4^k\, (n\vee M)^{k-1}\, ,
 $$

 \noindent where $H_n^{(k)}(\cdot)$ stands for the $k$-th derivative of $H_n(\cdot)$.

\end{theorem}

\medskip

\begin{corollary}
 \label{p:sXsub}

 Let $\gamma>0$. Let $(Z_n, \, n\ge 0)$ be a $\gamma$-dominable system.
 There exists a constant $c_5>0$, depending only on $(m, \, \gamma)$, such that for $M\ge M_0$, $n\ge 1$ and $v := m + \frac{1}{2c_4 (n\vee M)}$,
 \begin{equation}
     \E(Z_n^2 v^{Z_n})
     \le
     c_5 [\vartheta (n\vee M)]^{1/2} \prod_{i=0}^{n-1} [\E(m^{Z_i})]^{(m-1)/2} ;
     \label{cor1}
 \end{equation}
 in particular, we have, with $c_6 := \gamma^{1/2}c_5$,
 \begin{equation}
     \E(Z_n^2 v^{Z_n}) \le c_6 (n\vee M).
     \label{e:922745}
 \end{equation}

\end{corollary}

\medskip

The rest of the section is devoted to the proof of Theorem \ref{t:H_nkupper} and Corollary \ref{p:sXsub}. We start by mentioning a general technique going back to Collet et al.~\cite{collet-eckmann-glaser-martin}. Let $(X_n, \, n\ge 0)$ denote a Derrida--Retaux system satisfying $\E(X_0 \, m^{X_0})<\infty$. Let
$$
G_n(u) := \E(u^{X_n}) , \qquad n\ge 0.
$$

\noindent The iteration formula \eqref{iteration} is equivalent to:
\begin{equation}
    G_{n+1}(u) = \frac1u \, [G_n(u)]^m + (1-\frac1u) \, [G_n(0)]^m , \qquad n\ge 0.
    \label{iteration_Gn}
\end{equation}

\noindent A useful trick of Collet et al.~\cite{collet-eckmann-glaser-martin} consists in observing that this yields
$$
(u-1)u\, G_{n+1}'(u) - G_{n+1}(u)
=
[m(u-1)\, G_n'(u)-G_n(u)] \, [G_n(u)]^{m-1} \, .
$$

\noindent In particular, taking $u=m$ yields that
$$
G_{n+1}(m) - (m-1)m\, G_{n+1}'(m)
=
[G_n(m) - (m-1)m \, G_n'(m)] \, [G_n(m)]^{m-1} \, .
$$

\noindent Iterating this formula gives that for $n\ge 1$,
\begin{equation}
    G_n(m) - (m-1)m \, G_n'(m)
    =
    [G_0(m) - (m-1)m \, G_0'(m)] \prod_{i=0}^{n-1} [G_i(m)]^{m-1} \, .
    \label{iteration_Collet}
\end{equation}

\noindent In particular, \eqref{iteration_Collet} tells us that the sign of $\E(m^{X_n})-(m-1) \E(X_n\, m^{X_n})$ remains identical for all $n\ge 0$: it either is always positive (meaning that the system is supercritical), or is always negative (subcritical), or vanishes identically (critical).

A couple of known results which we are going to use for the subcritical or critical system: if $(X_n, \, n\ge 0)$ is a Derrida--Retaux satisfying $\E(m^{X_0}) \le (m-1) \E(X_0 \, m^{X_0}) < \infty$, then for all $n\ge 0$,
\begin{eqnarray}
    (m-1) \E(X_n \, m^{X_n})
 &\le& \E(m^{X_n})
    \le m^{1/(m-1)} ,
    \qquad n\ge 0,
    \label{Gn(m)<}
    \\
    \prod_{i=0}^{n-1} [\E(m^{X_i})]^{m-1}
 &\le& c_7 \, n^2,
    \qquad n\ge 1,
    \label{prod_Gi(m)<}
\end{eqnarray}

\noindent where $c_7>0$ is a constant depending on $m$ and on the law of $X_0$. See \cite[(3.11)]{bmvxyz_conjecture_DR} for the second inequality in \eqref{Gn(m)<} (proved in \cite{bmvxyz_conjecture_DR} for critical systems, and the same proof valid for subcritical systems as well), and \cite[Proposition 1]{bmxyz_questions} for \eqref{prod_Gi(m)<}.

The proof of Theorem \ref{t:H_nkupper} relies on the following preliminary result. Recall that $H_n(u) := \E(u^{Z_n})$.

\medskip

\begin{lemma}
\label{l:Pi_n<}

 Let $(Z_n, \, n\ge 0)$ be a $\gamma$-dominable system for some $\gamma>0$, in the sense of Definition $\ref{d:dominable}$.
 There exist constants $c_8>0$ and $c_9>0$, depending only on $m$, such that for $M\ge M_0$ and all integer $n\ge 1$,
 \begin{eqnarray}
     H_n''(m)
  &\le& c_8 \, [\vartheta (M)]^{1/2} \prod_{i=0}^{n-1} H_i(m)^{(m-1)/2},
     \label{Hn_2nd_derivative<}
     \\
     H_n'''(m)
  &\le& c_9 \, \vartheta (M) \prod_{i=0}^{n-1} H_i(m)^{m-1} .
     \label{Hn_3rd_derivative<}
 \end{eqnarray}

\end{lemma}

\medskip

\noindent {\it Proof.} Let
$$
D_n(u) := (m-1)u^3 H_n'''(u) + (4m-5)u^2 H_n''(u) + 2(m-2) uH_n'(u),
\qquad
n\ge 0.
$$

\noindent Recall from Remark \ref{r:dominable=>subcritical} that $m(m-1)H_0'(m) \le H_0(m)$. By \cite[Equation (19)]{bmxyz_questions}, this yields $D_{n+1}(m) \le D_n(m) H_n(m)^{m-1}$ for all $n\ge 0$. [In \cite{bmxyz_questions}, it was proved that if $m(m-1)H_0'(m) = H_0(m)$, then $D_{n+1}(m) = D_n(m) H_n(m)^{m-1}$, but the proof is valid for inequalities in place of equalities.] Accordingly, for all $n\ge 1$,
$$
D_n(m)
\le
D_0(m) \prod_{i=0}^{n-1} H_i(m)^{m-1} .
$$

\noindent By definition, $D_n(m) \ge m^3(m-1) H_n'''(m)$, so
\begin{equation}
    H_n'''(m) \le \frac{D_0(m)}{m^3(m-1)} \, \prod_{i=0}^{n-1} H_i(m)^{m-1}\, .
    \label{pf_Hn_2nd_derivative<_eq1}
\end{equation}

For any $j\ge 0$ and $\ell \in \{ 1, \, 2, \, 3\}$, writing $Z_j^3  m^{Z_j} \le \frac{9m^3}{2} \, Z_j(Z_j-1)(Z_j-2) m^{Z_j-3} \, {\bf 1}_{\{ Z_j \ge 3\}} + 8m^2$, and $Z_j^\ell \le Z_j^3$, we obtain
\begin{equation}
    \max_{\ell \in \{ 1, \, 2, \, 3\}} \E(Z_j^\ell  m^{Z_j}) \le \frac{9m^3}{2}\, H_j'''(m) + 8m^2\, .
    \label{pf_Hn_2nd_derivative<_eq2}
\end{equation}

\noindent In particular, with $j=0$ and $\ell \in \{ 1, \, 2\}$, this gives $\max\{ H_0'(m), \, H_0''(m)\} \le \frac{9m^3}{2}\, H_0'''(m) + 8m^2$. Consequently, with $c_{10} := m^3 (m-1)$, $c_{11} := m^2 (4m-5)$ and $c_{12} := 2m(m-2)$,
$$
D_0(m)
=
c_{10}\, H_0'''(m) + c_{11}\, H_0''(m) + c_{12}\, H_0'(m)
\le
c_{13} \, H_0'''(m) + c_{14}\, ,
$$

\noindent where $c_{13}:= c_{10} + \frac{9m^3}{2}(c_{11}+c_{12})$ and $c_{14}:=8m^2(c_{11}+c_{12})$. Since $H_0'''(m) \le \vartheta (M) +6 \le 7\vartheta (M)$ by assumption \eqref{assump_H0k} (applied to $k=3$; recalling that $\vartheta (M) \ge 1$), we get $D_0(m) \le 7c_{13} \, \vartheta (M) + c_{14} \le c_{15} \, \vartheta (M)$ with $c_{15} := 7c_{13} + c_{14}$. Going back to \eqref{pf_Hn_2nd_derivative<_eq1}, we have
$$
H_n'''(m) \le \frac{c_{15}}{m^3(m-1)} \, \vartheta (M) \prod_{i=0}^{n-1} H_i(m)^{m-1} ,
$$

\noindent proving \eqref{Hn_3rd_derivative<} with $c_9 := \frac{c_{15}}{m^3(m-1)}$.

It remains to prove \eqref{Hn_2nd_derivative<}. By \eqref{pf_Hn_2nd_derivative<_eq2} (applied to $j=n$),
$$
\E(Z_n^3 m^{Z_n}) \le \frac{9m^3}{2}\, H_n'''(m) + 8m^2\, ,
$$

\noindent whereas by \eqref{Gn(m)<}, $\E(Z_n m^{Z_n}) \le \frac{m^{1/(m-1)}}{m-1} =: c_{16}$, it follows from the Cauchy--Schwarz inequality that
$$
\E(Z_n^2 m^{Z_n}) \le c_{16}^{1/2} \Big( \frac{9m^3}{2}\, H_n'''(m) + 8m^2 \Big)^{\! 1/2}.
$$

\noindent Recall from \eqref{Hn_3rd_derivative<}, which we have just proved, that $H_n'''(m) \le c_9 \vartheta (M) \prod_{i=0}^{n-1} H_i(m)^{m-1}$. Writing $8m^2 \le 8m^2 \vartheta (M) \prod_{i=0}^{n-1} H_i(m)^{m-1}$ (because $\vartheta (M)\ge 1$ and $H_i(m) \ge 1$), this yields
$$
\E(Z_n^2 m^{Z_n})
\le
c_{17} \, [\vartheta (M)]^{1/2} \prod_{i=0}^{n-1} H_i(m)^{(m-1)/2},
$$

\noindent with $c_{17}:= c_{16}^{1/2} (\frac{9m^3}{2}\, c_9 + 8m^2)^{1/2}$. Since $H_n''(m) \le \E(Z_n^2 m^{Z_n})$, we obtain \eqref{Hn_2nd_derivative<} with $c_8 := c_{17}$.\qed

\medskip

\begin{remark}

 We often use the following inequalities for dominable systems:
 \begin{eqnarray}
     H_n''(m)
  &\le& c_8 \, [\vartheta (n\vee M)]^{1/2} \prod_{i=0}^{n-1} H_i(m)^{(m-1)/2},
     \label{Hn_2nd_derivative<bis}
     \\
     H_n'''(m)
  &\le& c_9 \, \vartheta (n\vee M) \prod_{i=0}^{n-1} H_i(m)^{m-1},
     \label{Hn_3rd_derivative<bis}
 \end{eqnarray}

 \noindent They are immediate consequences of Lemma \ref{l:Pi_n<} and the monotonicity of $M\mapsto \vartheta(M)$.\qed

\end{remark}

\medskip

We also need an elementary inequality.

\medskip

\begin{lemma}
 \label{l:201810171111}

 Let $\ell \ge 4$ be an integer, and let
 \begin{equation}
     B_\ell
     :=
     \{ {\bf u} := (u_1, \ldots, u_m) \in ([0, \, \ell-1] \cap \z)^m: \, u_1+\cdots+u_m = \ell \} \, .
     \label{B_ell}
 \end{equation}

 \noindent There exists a constant $c_{18}>0$, depending only on $m$, such that for all $y\ge 3m$,\footnote{Strictly speaking, we should write $\prod_{i: \, 1\le i\le m, \, u_i\ge 3} \frac{1}{u_i(u_i-1)}$ for $\prod_{i: \, u_i\ge 3} \frac{1}{u_i(u_i-1)}$. Notation: $\prod_\varnothing := 1$.}
 \begin{equation}
    \sum_{{\bf u} := (u_1, \ldots, u_m) \in B_\ell} y^{(\ell-\eta({\bf u})-2)^+} \prod_{i: \, u_i\ge 3} \frac{1}{u_i(u_i-1)}
    \le
    c_{18}\, \frac{y^{\ell-4}}{\ell^2} \, ,
    \label{pf_lemma51_eq2}
 \end{equation}

 \noindent where $a^+ := \max\{ a, \, 0\}$ as before, and
 \begin{equation}
     \eta({\bf u})
     :=
     \sum_{i=1}^m {\bf 1}_{\{ u_i \ge 1\} }
     \ge 2 \, .
     \label{eta}
 \end{equation}

\end{lemma}

\medskip

\noindent {\it Proof.} The sum over ${\bf u} := (u_1, \ldots, u_m) \in B_\ell$ satisfying $u_{\max} := \max_{1\le i\le m} u_i \le 2$ is very simple: in this case, $\ell \le 2m$; since $\eta({\bf u}) \ge 2$, we have $(\ell-\eta({\bf u})-2)^+ \le \ell-4$. The number of such ${\bf u}$ being smaller than $3^m$, we get
$$
\sum_{{\bf u}\in B_\ell: \, u_{\max} \le 2} y^{(\ell-\eta({\bf u}) -2)^+} \prod_{i: \, u_i\ge 3} \frac{1}{u_i(u_i-1)}
\le
3^m y^{\ell-4}
\le
c_{19} \, \frac{y^{\ell-4}}{\ell^2} \, ,
$$

\noindent with $c_{19} := 3^m (2m)^2$.  [Notation: $\sum_\varnothing := 0$.]

Let $\mathrm{LHS}_{\eqref{pf_lemma51_eq2}}$ denote the expression on the left-hand side of \eqref{pf_lemma51_eq2}. Then
$$
\mathrm{LHS}_{\eqref{pf_lemma51_eq2}}
\le
c_{19} \, \frac{y^{\ell-4}}{\ell^2}
+
\sum_{j=2}^{\ell \wedge m} \binom{m}{j} j!\, y^{(\ell-j-2)^+} \sum_{(u_1, \ldots, u_j)} \, \prod_{i: \, 1\le i\le j, \, u_i\ge 3} \frac{1}{u_i(u_i-1)} \, ,
$$

\noindent where, on the right-hand side, $\sum_{(u_1, \ldots, u_j)}$ sums over all $(u_1, \ldots, u_j) \in \z_+^j$ with $1\le u_1\le u_2\le \cdots\le u_j$ such that $u_1+\cdots+u_j=\ell$ and that $u_j \ge 3$. Note that $u_j\ge 3$ implies $j\le \ell-2$, thus $(\ell-j-2)^+ = \ell-j-2$. Moreover, we have $u_j \ge \frac{\ell}{j} \ge \frac{\ell}{m}$, thus $u_j(u_j-1) \ge \frac12 u_j^2 \ge \frac{\ell^2}{2m^2}$. Consequently, $\prod_{i\le j: \, u_i\ge 3} \frac{1}{u_i(u_i-1)}$ is bounded by $\frac{2m^2}{\ell^2} \prod_{i\le j-1: \, u_i\ge 3} \frac{1}{u_i(u_i-1)}$. This leads to (using $\binom{m}{j}j! \le m^j$):
\begin{eqnarray*}
    \mathrm{LHS}_{\eqref{pf_lemma51_eq2}}
 &\le& c_{19} \, \frac{y^{\ell-4}}{\ell^2}
    +
    \frac{2m^2}{\ell^2} \sum_{j=2}^{(\ell-2)\wedge m} m^j y^{\ell-j-2} \sum_{(u_1, \ldots, u_j)} \, \prod_{i\le j-1: \, u_i\ge 3} \frac{1}{u_i(u_i-1)}
    \\
 &\le& c_{19} \, \frac{y^{\ell-4}}{\ell^2}
    +
    \frac{2m^2}{\ell^2} \sum_{j=2}^{(\ell-2)\wedge m} m^j y^{\ell-j-2} \prod_{i=1}^{j-1} \Big( 1+\sum_{u=3}^\infty \frac{1}{u(u-1)} \Big) \, .
\end{eqnarray*}

\noindent Of course, $\sum_{u=3}^\infty \frac{1}{u(u-1)} = \frac12$; also, we bound $\sum_{j=2}^{(\ell-2)\wedge m}$ by $\sum_{j=2}^\infty$. This yields that
$$
\mathrm{LHS}_{\eqref{pf_lemma51_eq2}}
\le
c_{19} \, \frac{y^{\ell-4}}{\ell^2}
+
\frac{2m^2}{\ell^2} \sum_{j=2}^\infty m^j y^{\ell-j-2} (\frac32)^{j-1} \, .
$$

\noindent On the right-hand side, write $\sum_{j=2}^\infty m^j y^{\ell-j-2} (\frac32)^{j-1} = m^2 y^{\ell-4} \sum_{j=2}^\infty (\frac{m}{y})^{j-2} (\frac32)^{j-1}$; in view of our choice $y\ge 3m$, this is bounded by $m^2 y^{\ell-4} \sum_{j=2}^\infty (\frac13)^{j-2} (\frac32)^{j-1} = 3m^2 y^{\ell-4}$. As a consequence,
$$
\mathrm{LHS}_{\eqref{pf_lemma51_eq2}}
\le
c_{19} \, \frac{y^{\ell-4}}{\ell^2}
+
\frac{6m^4y^{\ell-4}}{\ell^2} \, ,
$$

\noindent yielding \eqref{pf_lemma51_eq2} with $c_{18} := c_{19}+6m^4$.\qed

\bigskip

We have all the ingredients for the proof of Theorem \ref{t:H_nkupper}.

\bigskip

\noindent {\it Proof of Theorem \ref{t:H_nkupper}.} Let $\gamma>0$ and let $(Z_n, \, n\ge 0)$ be a $\gamma$-dominable system. Write $H_n(u) := \E(u^{Z_n})$ as before. Recall $\vartheta(M) \ge 1$ (for all $M\ge M_0$) from \eqref{assump_H0k} and \eqref{assump_Qn}. Write, for brevity,
\begin{eqnarray*}
    Q_n
 &:=& \prod_{j=0}^{n-1} H_j(m)^{m-1}
    =
    \prod_{j=0}^{n-1} [\E(m^{Z_j})]^{m-1} ,
    \\
    M_n
 &:=& n \vee M ,
    \\
    \vartheta_n
 &:=& \vartheta(M_n)
    =
    \vartheta(n\vee M),
    \qquad
    n\ge 1, \; M\ge M_0.
\end{eqnarray*}

\noindent By assumption \eqref{assump_Qn},
\begin{equation}
    Q_n \le \vartheta_n Q_n \le \gamma \, M_n^2 ,
    \qquad
    n\ge 1 \, .
    \label{Qn}
\end{equation}

We claim that for all integer $k\ge 3$,
\begin{equation}
    H_n^{(k)}(m)
    \le
    c_4^{k-1}\, (k-2)! \, M_n^{k-3} \, \vartheta_n Q_n,
    \qquad
    n\ge 1,
    \label{e:momentupper1111}
\end{equation}

\noindent where $c_4 := \max\{ 4+\frac{c_{18}c_{20}}{m}, \, c_9^{1/2}, \, \gamma \}$, with $c_{20} := m^{m/(m-1)} (c_8^m \vee 1) (\gamma^{3m/2} \vee 1)$. Since $\vartheta_n Q_n \le \gamma \, M_n^2$ (see \eqref{Qn}), \eqref{e:momentupper1111} will imply Theorem \ref{t:H_nkupper}.

It remains to prove \eqref{e:momentupper1111}, which we do by induction in $k\ge 3$.

By \eqref{Hn_3rd_derivative<bis}, $H_n'''(m)\le c_9\, \vartheta_n Q_n$ for $n\ge 1$. So \eqref{e:momentupper1111} holds for $k=3$ since $c_9\le c_4^2$.

Let $\ell\ge 4$ be an integer. Suppose \eqref{e:momentupper1111} holds for all $k\in \{ 3, \, 4, \, \ldots, \, \ell-1\}$. We need to prove \eqref{e:momentupper1111} for $k=\ell$.

We first prove that the induction assumption yields that for $n\ge 1$ and ${\bf u} := (u_1, \ldots, u_m) \in B_\ell$ (defined in \eqref{B_ell}), we have, with $c_{20} := m^{m/(m-1)} (c_8^m \vee 1) (\gamma^{3m/2} \vee 1)$ as before,
\begin{equation}
    \prod_{i=1}^m H_n^{(u_i)}(m)
    \le
    c_{20} c_4^{\ell-2} M_n^{(\ell-\eta({\bf u}) -2)^+} \vartheta_n Q_n \prod_{i: \, u_i\ge 3} (u_i-2)! \, ,
    \label{e:diffmomentsupper}
\end{equation}

\noindent where $\eta({\bf u}):= \sum_{i=1}^m {\bf 1}_{\{ u_i \ge 1\} }$ is as in \eqref{eta}.

To check \eqref{e:diffmomentsupper}, let $n\ge 1$ and ${\bf u} \in B_\ell$. Since $H_n(m) \le m^{1/(m-1)} =: c_{21}$ (see \eqref{Gn(m)<}) and $H_n'(m) \le \frac{1}{m-1} \, \E(m^{Z_n}) \le \frac{c_{21}}{m-1} \le c_{21}$, we have
$$
\prod_{i=1}^m H_n^{(u_i)}(m)
\le
c_{21}^m H_n''(m)^{\lambda_2({\bf u})} \prod_{i: \, u_i\ge 3} H_n^{(u_i)}(m) \, ,
$$

\noindent where $\lambda_2({\bf u}) := \sum_{i=1}^m {\bf 1}_{\{ u_i =2\}}$. By \eqref{Hn_2nd_derivative<bis}, we have $H_n''(m) \le c_8 \, \vartheta_n^{1/2} \, Q_n^{1/2}$; thus with $c_{22} := c_{21}^m \max\{ c_8^m, \, 1\}$,
$$
\prod_{i=1}^m H_n^{(u_i)}(m)
\le
c_{22}\, (\vartheta_n Q_n)^{\lambda_2({\bf u})/2} \prod_{i: \, u_i\ge 3} H_n^{(u_i)}(m) \, .
$$

\noindent By the induction assumption in \eqref{e:momentupper1111}, $H_n^{(u_i)}(m) \le c_4^{u_i-1}(u_i-2)! \, M_n^{u_i-3}\, \vartheta_n Q_n$ if $u_i \ge 3$. As such, we have
\begin{eqnarray*}
    \prod_{i=1}^m H_n^{(u_i)}(m)
 &\le& c_{22}\, (\vartheta_n Q_n)^{\lambda_2({\bf u})/2} \prod_{i: \, u_i\ge 3} \Big( c_4^{u_i-1}(u_i-2)! \, M_n^{u_i-3} \vartheta_n Q_n\Big)
    \\
 &=& c_{22}\, (\vartheta_n Q_n)^{\lambda_2({\bf u})/2} \, (c_4 M_n)^{\sum_{i: \, u_i\ge 3} (u_i-1)} (M_n^{-2} \vartheta_n Q_n)^{\eta_3({\bf u})} \, \prod_{i: \, u_i\ge 3} (u_i-2)! \, ,
\end{eqnarray*}

\noindent where $\eta_3({\bf u}) := \sum_{i=1}^m {\bf 1}_{\{ u_i \ge 3\}}$. Note that $\sum_{i: \, u_i\ge 3} (u_i-1) = \sum_{i=1}^m (u_i-1)^+ - \lambda_2({\bf u}) = \ell - \eta({\bf u}) - \lambda_2({\bf u}) \le \ell-2$. So $c_4^{\sum_{i: \, u_i\ge 3} (u_i-1)} \le c_4^{\ell-2}$ (using $c_4>1$). This leads to:
\begin{eqnarray*}
    \prod_{i=1}^m H_n^{(u_i)}(m)
&\le& c_{22}\, c_4^{\ell-2}\, (\vartheta_n Q_n)^{\frac{\lambda_2({\bf u})}{2} + \eta_3({\bf u}) } \, M_n^{\ell - \eta({\bf u}) - \lambda_2({\bf u}) - 2\eta_3({\bf u})} \, \prod_{i: \, u_i\ge 3} (u_i-2)!
    \\
 &=& c_{22}\, c_4^{\ell-2}\, (\vartheta_n Q_n)^{\frac{\lambda_2({\bf u})}{2} + \eta_3({\bf u}) -1} \, M_n^{\ell - \eta({\bf u}) - \lambda_2({\bf u}) - 2\eta_3({\bf u})} \, \vartheta_n Q_n \prod_{i: \, u_i\ge 3} (u_i-2)! \; .
\end{eqnarray*}

Assume for the moment $\frac{\lambda_2({\bf u})}{2} + \eta_3({\bf u}) \ge 1$. By \eqref{Qn}, $\vartheta_n Q_n \le \gamma \, M_n^2$, so we have $(\vartheta_n Q_n)^{\frac{\lambda_2({\bf u})}{2} + \eta_3({\bf u}) -1} \le c_{23} \, M_n^{\lambda_2({\bf u}) + 2\eta_3({\bf u}) -2}$, with $c_{23} := \max\{ \gamma^{3m/2}, \, 1\}$. Since $c_{20} = c_{22}c_{23}$, we get
$$
\prod_{i=1}^m H_n^{(u_i)}(m)
\le
c_{20}\, c_4^{\ell-2}\, M_n^{\ell - \eta({\bf u}) - 2} \, \vartheta_n Q_n \prod_{i: \, u_i\ge 3} (u_i-2)! ,
$$

\noindent yielding \eqref{e:diffmomentsupper}. If, on the other hand, $\frac{\lambda_2({\bf u})}{2} + \eta_3({\bf u}) <1$, then $\lambda_2({\bf u}) \le 1$ and $\eta_3({\bf u}) =0$, i.e., $\max_{1\le i\le m} u_i\le 2$. This time, $\ell-\eta({\bf u}) = \sum_{i=1}^m (u_i-1)^+ \le 1$. The situation is very simple if we look at $\prod_{i=1}^m H_n^{(u_i)}(m)$ directly: at most one term among $H_n^{(u_i)}(m)$ is $H_n''(m)$ (which is bounded by $c_8 \, \vartheta_n^{1/2} \, Q_n^{1/2}$ as we have seen in \eqref{Hn_2nd_derivative<bis}), while all the rest is either $H_n(m)$ (which is bounded by $m^{1/(m-1)}=: c_{21}$) or $H_n'(m)$ (which is bounded by $1$ because by \eqref{Gn(m)<}, $H_n'(m) = \E(Z_n \, m^{Z_n-1}) \le \frac{m^{1/(m-1)}}{m(m-1)} \le 1$). Hence
$$
\prod_{i=1}^m H_n^{(u_i)}(m)
\le
c_{21}^m c_8 \, \vartheta_n^{1/2} \, Q_n^{1/2}
\le
c_{24} \, \vartheta_n Q_n \, ,
$$

\noindent with $c_{24} := c_{21}^m c_8$. [We have used $\vartheta_n Q_n \ge 1$.] This again gives \eqref{e:diffmomentsupper} because $c_{24} \le c_{20}$ and $c_4>1$.

Now that \eqref{e:diffmomentsupper} is proved, it is painless to complete the proof of Theorem \ref{t:H_nkupper}. Indeed, by \eqref{iteration_Gn},
$$
s H_{n+1}(s) = [H_n(s)]^m + (s-1) [H_n(0)]^m\, .
$$

\noindent On both sides, we differentiate $\ell$ times with respect to $s$ (recalling that $\ell \ge 4$, so the last term on the right-hand side, being affine in $s$, makes no contribution to the derivatives), and apply the general Leibniz rule to the first term on the right-hand side; this leads to:
\begin{eqnarray*}
    sH_{n+1}^{(\ell)}(s) + \ell H_{n+1}^{(\ell-1)}(s)
 &=& \sum_{(u_1, \ldots, u_m) \in \z_+^m: \, u_1+\cdots+u_m=\ell} \frac{\ell!}{u_1! \cdots u_m!} \prod_{i=1}^m H_n^{(u_i)}(s)
    \\
 &=& m H_n^{(\ell)}(s) H_n(s)^{m-1}
    +
    \sum_{{\bf u}\in B_\ell} \frac{\ell!}{u_1! \cdots u_m!} \prod_{i=1}^m H_n^{(u_i)}(s)\, .
\end{eqnarray*}

\noindent Note that the expression on the left-hand side is at least $sH_{n+1}^{(\ell)}(s)$. We take $s=m$ to see that
$$
H_{n+1}^{(\ell)}(m)
-
H_n^{(\ell)}(m) H_n(m)^{m-1}
\le
\frac1m \sum_{{\bf u}\in B_\ell} \frac{\ell!}{u_1! \cdots u_m!} \prod_{i=1}^m H_n^{(u_i)} (m).
$$

\noindent By \eqref{e:diffmomentsupper}, we have $\prod_{i=1}^m H_n^{(u_i)}(m) \le c_{20} c_4^{\ell-2} M_n^{(\ell-\eta({\bf u})-2)^+} \vartheta_n Q_n \prod_{i: \, u_i\ge 3} (u_i-2)!$, where $\eta ({\bf u}) := \sum_{i=1}^m {\bf 1}_{\{ u_i \ge 1\} }$ is as in \eqref{eta}. Hence
$$
H_{n+1}^{(\ell)}(m)
-
H_n^{(\ell)}(m) H_n(m)^{m-1}
\le
\frac{c_{20}}{m} c_4^{\ell-2} \ell! \, \vartheta_n Q_n \sum_{{\bf u}\in B_\ell} M_n^{(\ell-\eta({\bf u})-2)^+} \prod_{i: \, u_i\ge 3} \frac{1}{u_i(u_i-1)} \, ,
$$

\noindent which, in view of Lemma \ref{l:201810171111} (applied to $y:= M_n$), yields that, for $n\ge 1$,
$$
H_{n+1}^{(\ell)}(m)
\le
H_n^{(\ell)}(m)H_n(m)^{m-1}
+
c_{18} \frac{c_{20}}{m} c_4^{\ell-2}\, (\ell-2)! \, M_n^{\ell-4} \vartheta_n Q_n.
$$

\noindent Recall that $Q_n := \prod_{j=0}^{n-1} H_j(m)^{m-1}$. Iterating this inequality, and by means of the monotonicity of $n\mapsto M_n^{\ell-4} \vartheta_n Q_n$, we get
\begin{eqnarray*}
    H_n^{(\ell)}(m)
 &\le& H_0^{(\ell)}(m) Q_n
    +
    \sum_{j=0}^{n-1} \frac{c_{18}c_{20}}{m} c_4^{\ell-2} (\ell-2)! \, M_n^{\ell-4} \, \vartheta_n Q_n
    \\
 &=& \Big( \frac{H_0^{(\ell)}(m)}{\vartheta_n} +n\frac{c_{18}c_{20}}{m} c_4^{\ell-2} (\ell-2)! \, M_n^{\ell-4} \Big) \vartheta_n Q_n \, .
\end{eqnarray*}

\noindent We use $n \le M_n$ so that $n\, M_n^{\ell-4} \le M_n^{\ell-3}$. On the other hand, $H_0^{(\ell)}(m) \le M^{\ell-3} (\vartheta(M) + \ell !)$ (by assumption \eqref{assump_H0k}), which is bounded by $M_n^{\ell-3} \vartheta_n\, (1+ \ell !)$. Thus
\begin{eqnarray*}
    H_n^{(\ell)}(m)
 &\le& \Big( 1+ \ell! + \frac{c_{18}c_{20}}{m} c_4^{\ell-2} (\ell-2)! \Big) M_n^{\ell-3} \, \vartheta_n Q_n
    \\
 &\le& (1 + \ell(\ell-1) +\frac{c_{18}c_{20}}{m} c_4^{\ell-2}) \, (\ell-2)! \, M_n^{\ell-3} \, \vartheta_n Q_n \, .
\end{eqnarray*}

\noindent Since $\ell \ge 4$, we have $1 + \ell(\ell-1) \le \ell^2 \le 2^\ell \le 4 \, c_4^{\ell-2}$ (because $c_4\ge 2$), so $1 + \ell(\ell-1) +\frac{c_{18}c_{20}}{m} c_4^{\ell-2} \le 4 \, c_4^{\ell-2} + \frac{c_{18}c_{20}}{m} c_4^{\ell-2}\le c_4^{\ell-1}$ by means of the fact that $c_4 := \max\{ 4+\frac{c_{18}c_{20}}{m}, \, c_9^{1/2}, \, \gamma \}$. Consequently,
$$
H_n^{(\ell)}(m)
\le
c_4^{\ell-1} (\ell-2)! \, M_n^{\ell-3} \, \vartheta_n Q_n \, ,
$$

\noindent implying \eqref{e:momentupper1111} for $k=\ell$, and completing the proof of Theorem \ref{t:H_nkupper}.\qed

\bigskip

\noindent {\it Proof of Corollary \ref{p:sXsub}.} Only \eqref{cor1} needs proving because \eqref{e:922745} will follow immediately from \eqref{cor1} and assumption \eqref{assump_Qn}.

Let $n\ge 1$ and $s\in [m, \, m + \frac{1}{2c_4 M_n}]$, where we keep using the notation
$$
M_n := n \vee M\, .
$$

\noindent Write $H_n(u) := \E(u^{Z_n})$ as before. Then
$$
\E(Z_n^2 s^{Z_n})
=
s H_n'(s) + s^2H_n''(s) \, .
$$

\noindent Since $u\mapsto H_n''(u)$ is non-decreasing, we have $H_n'(s) \le H_n'(m) + (s-m) H_n''(s)$; hence
$$
\E(Z_n^2 s^{Z_n})
\le
sH_n'(m) + s(s-m)H_n''(s) + s^2H_n''(s)
=
sH_n'(m) + s(2s-m)H_n''(s)\, .
$$

\noindent On the right-hand side, we use $H_n'(m) \le 1$ (which has already been observed as a consequence of \eqref{Gn(m)<}), and $m\le s\le m+1$ (so $s(2s-m) \le (m+1)(m+2)$), to see that
$$
\E(Z_n^2 s^{Z_n})
\le
m+1 + (m+1)(m+2) H_n''(m+\frac{1}{2c_4 M_n}) \, .
$$

\noindent To bound $H_n''(m+\frac{1}{2c_4 M_n})$, we recall that $Z_n$ is bounded for each $n$ (which is a consequence of the boundedness of $Z_0$), so by Taylor expansion,
$$
H_n''(m+\frac{1}{2c_4 M_n})
-
H_n''(m)
=
\sum_{k=3}^\infty \frac{(2c_4 M_n)^{-(k-2)}}{(k-2)!} H_n^{(k)}(m) \, .
$$

\noindent By \eqref{e:momentupper1111}, we have $H_n^{(k)}(m)\le c_{4}^{k-1}\, (k-2)! \, M_n^{k-3} \, \vartheta_n Q_n$ (for $k\ge 3$). Hence
$$
H_n''(m+\frac{1}{2c_4 M_n})
-
H_n''(m)
\le
\sum_{k=3}^\infty \frac{(2c_4 M_n)^{-(k-2)}}{(k-2)!} c_{4}^{k-1}\, (k-2)! \, M_n^{k-3} \, \vartheta_n Q_n
=
c_4\, \frac{\vartheta_n Q_n}{M_n} \, .
$$

\noindent As such, we arrive at:
$$
\E(Z_n^2 s^{Z_n})
\le
m+1 + (m+1)(m+2) \Big( H_n''(m) + c_4\, \frac{\vartheta_n Q_n}{M_n} \Big) \, .
$$

\noindent By \eqref{Hn_2nd_derivative<bis}, $H_n''(m)\le c_8 \, \vartheta_n^{1/2} Q_n^{1/2}$, whereas according to assumption \eqref{assump_Qn}, $M_n \ge \frac{1}{\gamma^{1/2}} \, \vartheta_n^{1/2} Q_n^{1/2}$, this readily yields \eqref{cor1} (recalling that $\vartheta_n Q_n \ge 1$). Corollary \ref{p:sXsub} is proved.\qed

\section{Proof of Theorem \ref{t:moment_ub}}
\label{s:moment_s>m}

Let $(Y_n, \, n\ge 0)$ be a critical system such that $\E(s^{Y_0})<\infty$ for some $s>m$. We claim that with $c_{25} := (\sup_{x>0} x\, (\frac{s}{m})^{-x/\ee}) \vee 1 \in [1, \, \infty)$, we have, for all integer $k\ge 1$,
\begin{equation}
    \E(Y_0^k m^{Y_0}) \le c_{25}^k \, \E(s^{Y_0}) \, k!\, .
    \label{case_big_moments}
\end{equation}

Indeed, by definition,
$$
x\le c_{25} \, (\frac{s}{m})^{x/\ee} \, ,
$$

\noindent for all $x>0$. Taking to the power $k$ on both sides and with $x:= \frac{\ee\, \ell}{k}$, we see that for all integers $k\ge 1$ and $\ell \ge 1$,
$$
(\frac{\ee\, \ell}{k})^k \le c_{25}^k \, (\frac{s}{m})^\ell \, .
$$

\noindent Since $k! \ge (\frac{k}{\ee})^k$ by Stirling's formula, this yields $\ell^k m^\ell\le c_{25}^k \, s^\ell \, k!$, from which \eqref{case_big_moments} follows.

Let $L\ge 1$ be an integer, and let $Z_0 = Z_0(M, \, L) := Y_0 \, {\bf 1}_{\{ Y_0 \le L\} }$. Then $Z_0$ is bounded, and does not depend on $M$, though we still treat it as indexed by $M$. Let us check that assumptions \eqref{assump_H0k} and \eqref{assump_Qn} in Definition \ref{d:dominable} are satisfied.

Since $Z_0\le Y_0$, it follows from \eqref{case_big_moments} that $\E(Z_0^k m^{Z_0}) \le c_{26} c_{25}^k \, k!$ for $k\ge 1$, where $c_{26} =c_{26}(s) := \E(s^{Y_0}) \ge 1$. Let $M\ge c_{26} c_{25}^4$. Then $M^{k-3} \ge (c_{26} c_{25}^3)^{k-3} c_{25}^{k-3} \ge c_{26} c_{25}^3 c_{25}^{k-3} = c_{26} c_{25}^k$ for all integer $k\ge 4$, so $c_{26} c_{25}^k \, k! \le M^{k-3} k!$ for $k\ge 4$. For $k=3$, we have $c_{26} c_{25}^3 \, 3! \le c_{27} + 3!$ with $c_{27} := 6(c_{26} c_{25}^3-1) \vee 1 \in [1, \, \infty)$. As such, we see that for all integers $k\ge 3$ and $M\ge c_{26} c_{25}^4$,
$$
\E(Z_0^k m^{Z_0}) \le M^{k-3} (c_{27} + k!) ;
$$

\noindent in words, assumption \eqref{assump_H0k} is satisfied with $\vartheta (M) := c_{27}$ for all $M \ge M_0 := \lceil c_{26} c_{25}^4\rceil$.

Assumption \eqref{assump_Qn} is easily seen to be satisfied: by \eqref{prod_Gi(m)<}, $\prod_{i=0}^{n-1} [\E(m^{Y_i})]^{m-1} \le c_7 n^2$ for all $n\ge 1$, so \eqref{assump_Qn} holds with $\gamma := c_7 c_{27}$.

So we are entitled to apply Theorem \ref{t:H_nkupper} to see that for $k\ge 3$, $n\ge 1$ and $M\ge M_0:= \lceil c_{26} c_{25}^4\rceil$,
$$
\E[Z_n(Z_n-1) \cdots (Z_n-k+1)\, m^{Z_n-k}] \le k! \, c_4^k \, (n\vee M)^{k-1} \, .
$$

\noindent We take $M := M_0$. Recall that $Z_0 = Z_0(M_0, \, L) := Y_0 \, {\bf 1}_{\{ Y_0 \le L\} }$. Letting $L\to \infty$, and applying the monotone convergence theorem, we get, for $k\ge 3$ and $n\ge 1$,
$$
\E[Y_n(Y_n-1) \cdots (Y_n-k+1)\, m^{Y_n-k}] \le k! \, c_4^k \, (n\vee M_0)^{k-1} \, .
$$

\noindent This implies the desired inequality for $k\ge 3$. The case $k=1$ has already been implicitly treated in the proof of Corollary \ref{p:sXsub} in Section \ref{s:domination}: $\E(Y_n \, m^{Y_n-1}) = \frac{\E(m^{Y_n})}{m(m-1)} \le \frac{m^{1/(m-1)}}{m(m-1)}$ (see \eqref{Gn(m)<}). The case $k=2$ follows from the cases $k=1$ and $k=3$ by the Cauchy--Schwarz inequality.\qed

\section{Truncating the critical system}
\label{s:truncating}

Let $(Y_n, \, n\ge 0)$ be a Derrida--Retaux system satisfying $\E(m^{Y_0}) = (m-1) \E(Y_0\, m^{Y_0})<\infty$ (so the system is critical). Recall that the system $(Y_n, \, n\ge 0)$ is of ``finite variance" if $\E(Y_0^3 \, m^{Y_0})<\infty$, and is stable if $\P(Y_0=j) \sim c_0 \, m^{-j}j^{-\alpha}$, $j\to \infty$, for some $c_0>0$ and $2<\alpha<4$ as in \eqref{alpha}. The following theorem tells that in either case, we can define
$$
Z_0 := Y_0 \, {\bf 1}_{\{ Y_0 \le a(M)\}},
$$

\noindent for some appropriate $a(M) \in [1, \, \infty]$ such that $(Z_n, \, n\ge 0)$ is dominable in the sense of Definition \ref{d:dominable}; the values of $a(M)$ and $\vartheta(M)$ (as defined in Definition \ref{d:dominable}) are also given as they are often useful in the applications.

\medskip

\begin{theorem}
\label{t:dominable}

 Let $(Y_n, \, n\ge 0)$ be a Derrida--Retaux system satisfying $\E(m^{Y_0}) = (m-1) \E(Y_0\, m^{Y_0})<\infty$. If it is either of ``finite variance" or stable, then there is a dominable system $(Z_n, \, n\ge 0)$ such that $Z_0 \le Y_0$ \hbox{\rm a.s.} More precisely,

 {\rm (i)} if the system is of ``finite variance", we can choose $Z_0 := Y_0 \, {\bf 1}_{\{ Y_0 \le M \zeta(M)\} }$ where $\zeta(M) := - \log \E(Y_0^3 m^{Y_0} \, {\bf 1}_{\{ Y_0 >M\} }) \le \infty$,\footnote{So $\lim_{M\to \infty} \zeta(M) =\infty$ by the ``finite variance" assumption.} with $\vartheta (M) := \max\{ \E(Y_0^3 m^{Y_0}), \, 1\}$;

 {\rm (ii)} if the system is stable, we can choose $Z_0 := Y_0 \, {\bf 1}_{\{ Y_0 \le M\} }$, with $\vartheta(M) := c_{30}\, M^{4-\alpha}$, where $c_{30}$ is the constant in \eqref{theta(M)_stable} below.

\end{theorem}

\medskip

For the sake of clarity, the two situations (``finite variance", stable) are discussed in distinct parts.

\subsection{Proof of Theorem \ref{t:dominable}: the ``finite variance" case}

Assume $\E(m^{Y_0}) = (m-1) \E(Y_0\, m^{Y_0})<\infty$ and $\E(Y_0^3 \, m^{Y_0}) <\infty$. Let
$$
\zeta(M)
:=
- \log \E(Y_0^3 m^{Y_0} \, {\bf 1}_{\{ Y_0 >M\} })
\le
\infty \, .
$$

\noindent Since $\lim_{M\to \infty} \zeta(M) \to \infty$, we can choose $M$ sufficiently large so that $M \, \zeta(M) \ge 2$. Let $Z_0 := Y_0 \, {\bf 1}_{\{ Y_0 \le M \zeta(M)\} }$.

We claim that assumption \eqref{assump_H0k} is satisfied with $\vartheta (M) := \max\{ \E(Y_0^3 m^{Y_0}), \, 1\}$ and that $Z_0$ is bounded. For any integer $k\ge 3$, we write
\begin{eqnarray*}
   \E(Z_0^km^{Z_0})
 &=& \E(Y_0^k m^{Y_0} \, {\bf 1}_{\{Y_0\le M\}})
   +
   \E(Y_0^k m^{Y_0} \, {\bf 1}_{\{M< Y_0\le M\zeta(M)\}})
   \\
 &\le& M^{k-3}\E(Y_0^3 m^{Y_0} \, {\bf 1}_{\{Y_0\le M\}})
   +
   \E(Y_0^k m^{Y_0} \, {\bf 1}_{\{M< Y_0\le M\zeta(M)\}}) \, .
\end{eqnarray*}

\noindent The first term on the right-hand side is easy to handle: we have $\E(Y_0^k m^{Y_0} \, {\bf 1}_{\{Y_0\le M\}}) \le M^{k-3}\E(Y_0^3 m^{Y_0} \, {\bf 1}_{\{Y_0\le M\}}) \le \E(Y_0^3 m^{Y_0}) \, M^{k-3}$. In case $\zeta(M) =\infty$, we have $Y_0\le M$ a.s., so $Z_0$ is bounded and the second term on the right-hand side vanishes, which yields \eqref{assump_H0k} with $\vartheta (M) := \max\{ \E(Y_0^3 m^{Y_0}), \, 1\}$.

To treat the case $\zeta(M)<\infty$ (in which case $Z_0$ is obviously bounded), let us look at the second term on the right-hand side: since $\E(Y_0^3 m^{Y_0} \, {\bf 1}_{\{ Y_0>M\}}) = \ee^{-\zeta(M)}$, we have
$$
\E(Y_0^k m^{Y_0} \, {\bf 1}_{\{M< Y_0\le M\zeta(M)\}})
\le
M^{k-3}\zeta(M)^{k-3}\E(Y_0^3 m^{Y_0} \, {\bf 1}_{\{ Y_0>M\}})
=
M^{k-3}\zeta(M)^{k-3}\, \ee^{-\zeta(M)}\, .
$$

\noindent Applying the inequality $\ee^x \ge \frac{x^{k-3}}{(k-3)!}$ (for $x\ge 0$ and $k\ge 3$) to $x:= \zeta(M)$ yields \eqref{assump_H0k} again with $\vartheta (M) := \max\{ \E(Y_0^3 m^{Y_0}), \, 1\}$.

Consequently, regardless of whether $\zeta(M)$ is finite or infinite, $Z_0$ is bounded, and assumption \eqref{assump_H0k} is satisfied with $\vartheta (M) := \max\{ \E(Y_0^3 m^{Y_0}), \, 1\}$. Note that $\vartheta(M)$ does not depend on $M$.

Assumption \eqref{assump_Qn} is also satisfied: by \eqref{prod_Gi(m)<}, $\prod_{i=0}^{n-1} [\E(m^{Y_i})]^{m-1} \le c_7 n^2$ for all $n\ge 1$; since $\E(m^{Z_i}) \le \E(m^{Y_i})$ for all $i\ge 0$, \eqref{assump_Qn} is satisfied with $\gamma := c_7 \max\{ \E(Y_0^3 m^{Y_0}), \, 1\}$.\qed

\subsection{Proof of Theorem \ref{t:dominable}: the stable case}

We start with a simple inequality.

\medskip

\begin{lemma}
\label{l:prod_subcritical}

 Let $(X_n, \, n\ge 0)$ be a Derrida--Retaux system satisfying $(m-1)\E(X_0 \, m^{X_0}) <\E(m^{X_0})  < \infty$. Then
 $$
 \prod_{i=0}^\infty [\E(m^{X_i})]^{m-1}
 \le
 \frac{1}{\E(m^{X_0})- (m-1) \E(X_0\, m^{X_0})} \, .
 $$

\end{lemma}

\medskip

\noindent {\it Proof.} For the moment, let $(X_n, \, n\ge 0)$ be an arbitrary Derrida--Retaux system satisfying $\E(X_0 \, m^{X_0})<\infty$, and such that $\E(m^{X_0}) \not= (m-1) \E(X_0\, m^{X_0})$. By \eqref{iteration_Collet},
\begin{equation}
   \prod_{i=0}^{n-1} [\E(m^{X_i})]^{m-1}
    =
   \frac{\E(m^{X_n})- (m-1)\E(X_n\, m^{X_n})}{\E(m^{X_0})- (m-1) \E(X_0\, m^{X_0})} \, .
   \label{recursion_Coll}
\end{equation}

\noindent For the nominator in \eqref{recursion_Coll}, we observe that $\E(m^{X_n})- (m-1)\E(X_n\, m^{X_n}) = \E[(1-(m-1)X_n) m^{X_n}] \le \P(X_n=0) \le 1$. If we assume $(m-1)\E(X_0 \, m^{X_0}) <\E(m^{X_0})$, then the denominator in \eqref{recursion_Coll} is positive; the lemma follows immediately from the monotone convergence theorem by letting $n\to \infty$.\qed

\bigskip

We now proceed to the proof of Theorem \ref{t:dominable} for stable systems. Assume $\E(m^{Y_0}) = (m-1) \E(Y_0\, m^{Y_0})<\infty$ and $\P(Y_0=j) \sim c_0 \, m^{-j}j^{-\alpha}$, $j\to \infty$, for some $c_0>0$ and $2<\alpha<4$ as in \eqref{alpha}. This yields the existence of constants $c_{28} \ge c_{29} >0$ and an integer $j_0\ge 1$, all depending on $m$ and on the law of $Y_0$, such that
\begin{eqnarray}
    \P(Y_0=j)
 &\le& c_{28} \, m^{-j}j^{-\alpha} ,
    \qquad
    j\ge 1,
    \label{stable_tail_ub}
    \\
    \P(Y_0=j)
 &\ge& c_{29} \, m^{-j}j^{-\alpha} .
    \qquad
    j\ge j_0,
    \label{stable_tail_lb}
\end{eqnarray}

\noindent Let $M\ge j_0$ be an integer and let
\begin{equation}
    Z_0 := Y_0 \, {\bf 1}_{\{ Y_0 \le M\} } ,
    \label{Z0}
\end{equation}

\noindent which is a bounded random variable. For integer $k\ge 3$, we write
$$
\E(Z_0^k m^{Z_0})
=
\E(Y_0^k m^{Y_0} \, {\bf 1}_{\{Y_0\le M\}})
\le
M^{k-3} \, \E(Y_0^3 m^{Y_0} \, {\bf 1}_{\{Y_0\le M\}}) ,
$$

\noindent so by \eqref{stable_tail_ub}, we have
\begin{eqnarray*}
    \E(Z_0^k m^{Z_0})
 &\le& c_{28} M^{k-3} \sum_{j=1}^M j^{3-\alpha}
    \le
    c_{28} M^{k-3} \int_0^{M+1} x^{3-\alpha} \d x
    \\
 &=& c_{28} M^{k-3} \, \frac{(M+1)^{4-\alpha}}{4-\alpha}
    \le
    c_{28} M^{k-3} \, \frac{2^{4-\alpha} \, M^{4-\alpha}}{4-\alpha}  \, .
\end{eqnarray*}

\noindent As such, assumption \eqref{assump_H0k} is satisfied with
\begin{equation}
    \vartheta (M) := c_{30}\, M^{4-\alpha},
    \label{theta(M)_stable}
\end{equation}

\noindent where $c_{30} := \max\{ c_{28} \frac{2^{4-\alpha}}{4-\alpha}, \, 1\}$. In particular, $M\mapsto \vartheta (M)$ is non-decreasing.

It remains to check assumption \eqref{assump_Qn}, which in this case states that for some constant $c_{31}>0$,
\begin{equation}
    \prod_{i=0}^{n-1} [\E(m^{Z_i})]^{m-1}
    \le
    c_{31} \, (n\vee M)^{\alpha-2},
    \qquad
    M\ge j_0, \; n\ge 1.
    \label{check_assump2_stable}
\end{equation}

Since $\E(m^{Z_i}) \le \E(m^{Y_i})$, there is nothing to prove if $n \le j_0$: it suffices to take $c_{31}$ such that $c_{31} \ge \prod_{i=0}^{j_0-1} [\E(m^{Y_i})]^{m-1}$. Let us assume $n>j_0$. We write $\prod_{i=0}^{n-1} [\E(m^{Z_i})]^{m-1} \le \prod_{i=0}^\infty [\E(m^{Z_i})]^{m-1}$, so by Lemma \ref{l:prod_subcritical},
$$
\prod_{i=0}^{n-1} [\E(m^{Z_i})]^{m-1}
\le
\frac{1}{\E(m^{Z_0})- (m-1) \E(Z_0\, m^{Z_0})} \, .
$$

\noindent By definition, $Z_0 = Y_0 \, {\bf 1}_{\{ Y_0 \le M\} }$, and by assumption, $\E(m^{Y_0}) = (m-1) \E(Y_0\, m^{Y_0})$. So
\begin{eqnarray}
 &&\E(m^{Z_0})- (m-1) \E(Z_0\, m^{Z_0})
    \nonumber
    \\
 &=& \E[((m-1)Y_0 -1)\, m^{Y_0}\, {\bf 1}_{\{ Y_0 > M\} }) + \P(Y_0 > M)
    \label{E(mZ0)}
    \\
 &\ge& \E[((m-1)Y_0 -1)\, m^{Y_0}\, {\bf 1}_{\{ Y_0 > M\} }) \, ,
    \nonumber
\end{eqnarray}

\noindent which, by \eqref{stable_tail_lb}, is $\ge c_{29} \sum_{j=M+1}^\infty ((m-1)j -1) j^{-\alpha} \ge \frac{c_{32}}{M^{\alpha-2}}$ for some constant $c_{32}>0$ and all $M\ge j_0$. Consequently,
\begin{equation}
    \prod_{i=0}^{n-1} [\E(m^{Z_i})]^{m-1}
    \le
    \frac{M^{\alpha-2}}{c_{32}} \, ,
    \label{e:Q_nupper_alpha}
\end{equation}

\noindent which is bounded by $\frac{(n\vee M)^{\alpha-2}}{c_{32}}$ as $\alpha>2$. This yields \eqref{check_assump2_stable}.\qed

\medskip

\begin{remark}
\label{r:Z2upperQ}

 For further use, let us note that in the stable case, $\vartheta (M) := c_{30}\, M^{4-\alpha}$ for all $M\ge M_0$ (by \eqref{theta(M)_stable}), whereas $\prod_{i=0}^{n-1} [\E(m^{Z_i})]^{m-1} \le \frac{M^{\alpha-2}}{c_{32}}$ (by \eqref{e:Q_nupper_alpha}), so inequality \eqref{cor1} in Corollary \ref{p:sXsub} implies that for $M\ge M_0$, $n\ge 1$ and $v := m + \frac{1}{2c_4 (n\vee M)}$ ($c_4 \ge 1$ being as before the constant in Theorem $\ref{t:H_nkupper}$),
\begin{equation}
    \E(Z_n^2 v^{Z_n})
    \le
    c_{33} (n\vee M)^{(4-\alpha)/2} M^{(\alpha-2)/2} \, ,
    \label{cor1_variant_stable}
\end{equation}

\noindent where $c_{33} := \frac{c_5 \, c_{30}^{1/2}}{c_{32}^{1/2}}$.\qed

\end{remark}

\section{Proof of Theorem \ref{t:main}}
\label{s:pf_thm}

Let $(Y_n, \, n\ge 0)$ be a Derrida--Retaux system such that $\E(m^{Y_0}) =(m-1) \,\E(Y_0 \, m^{Y_0})<\infty$. We assume that $\P(Y_0=j) \sim c_0 \, m^{-j}j^{-\alpha}$, $j\to \infty$, for some $c_0>0$ and $2<\alpha<4$ as in \eqref{alpha}.

\subsection{Upper bound}

We start with a lemma.

\medskip

\begin{lemma}
\label{l:general<}

 Let $(X_n, \, n\ge 0)$ be a Derrida--Retaux system satisfying $\E(X_0 \, m^{X_0}) < \infty$. There exists a constant $c_{34} >0$, depending only on $m$, such that $(\frac{n^2}{0} := \infty)$
 $$
 \prod_{i=0}^{n-1} [\E(m^{X_i})]^{m-1}
 \le
 c_{34} \, \frac{n^2}{\E(X_0^3 m^{X_0} {\bf 1}_{\{ 2\le X_0 \le 3n\} })} ,
 \qquad
 n\ge 1 \, .
 $$

\end{lemma}

\medskip

\noindent {\it Proof.} The lemma was known in various forms. Recall from \cite{bmxyz_questions} that for $s\in (\frac{m}{2}, \, m)$,
\begin{equation}
    \prod_{i=0}^{n-1} [\E(m^{X_i})]^{m-1}
    \le
    \Big(\frac{m}{2s-m}\Big)^n \frac{1}{\Delta_0(s)},
    \label{tech1}
\end{equation}

\noindent where
$$
\Delta_0(s)
:=
\sum_{k=1}^\infty m^k ((m-1)k-1) (1-(k+1)x^k + kx^{k+1}) \P(X_0=k),
$$

\noindent with $x:= \frac{s}{m}\in (\frac12, \, 1)$. We now reproduce some elementary computations from \cite{bmvxyz_conjecture_DR}. For $k\ge 1$, we have $x^k \le \ee^{-(1-x)k}$, so $1-(1+k)x^k+ kx^{1+k} \ge 1-(1+u)\, \ee^{-u}$, where $u:= (1-x)k>0$. Since $1 - (1+v)\,\ee^{-v} \ge 1 -\frac{2}{\ee}$ for $v\ge 1$ (because $v\mapsto 1 - (1+v)\,\ee^{-v}$ is increasing on $(0, \, \infty)$) and $1 - (1+v)\, \ee^{-v} \ge \frac{v^2}{2 \ee}$ for $v\in (0,\, 1]$ (because $v \mapsto 1 - (1+v)\, \ee^{-v} - \frac{v^2}{2\ee}$ is increasing on $(0, \, 1]$), we get, for $k\ge 1$,
$$
1-(k+1)x^k + kx^{k+1}
\ge
c_{35} \, \min\{ (1-x)^2 k^2, \, 1\} \, ,
$$

\noindent where $c_{35} := \min\{ 1 -\frac{2}{\ee}, \, \frac{1}{2\ee}\} > 0$. Let $n\ge 1$. We take $s=s_n:= (1-\frac{1}{3n})m$ (so $x= 1-\frac{1}{3n}$), to see that
\begin{eqnarray*}
    \Delta_0(s_n)
 &\ge& c_{35} \sum_{k=1}^\infty m^k ((m-1)k-1) \min\{ \frac{k^2}{(3n)^2}, \, 1\} \P(X_0=k)
    \\
 &\ge& \frac{c_{35}}{(3n)^2} \sum_{k=1}^{3n} k^2 m^k ((m-1)k-1) \P(X_0=k)\, .
\end{eqnarray*}

\noindent We use $(m-1)k-1 \ge \frac{m-1}{2} \, k$ for $k\ge 2$, so that with $c_{36} := c_{35} \, \frac{m-1}{18} >0$
$$
\Delta_0(s_n)
\ge
\frac{c_{36}}{n^2} \sum_{k=2}^{3n} k^3 m^k \P(X_0=k)
=
\frac{c_{36}}{n^2} \, \E(X_0^3 m^{X_0} \, {\bf 1}_{\{ 2\le X_0 \le 3n\} } ).
$$

\noindent This, in view of \eqref{tech1} (applied to $s=s_n$), yields the lemma.\qed

\bigskip

\noindent {\it Proof of the upper bound in Theorem \ref{t:main}.} By Lemma \ref{l:general<}, for all $n\ge 1$,
$$
\prod_{i=0}^{n-1} [\E(m^{Y_i})]^{m-1}
\le
c_{34} \, \frac{n^2}{\E(Y_0^3 m^{X_0} {\bf 1}_{\{ 2\le Y_0 \le 3n\} })} .
$$

\noindent Since $\P(Y_0=j) \sim c_0 \, m^{-j}j^{-\alpha}$, $j\to \infty$, and $2<\alpha<4$, we have $\E(Y_0^3 m^{Y_0} {\bf 1}_{\{ 2\le Y_0 \le 3n\} }) \ge c_{37}\, n^{4-\alpha}$ for some constant $c_{37}>0$ and all sufficiently large $n$; this implies the upper bound in Theorem \ref{t:main}.\qed

\subsection{Lower bound in Theorem \ref{t:main}}

The proof of the lower bound in Theorem \ref{t:main} needs some preparation.

\medskip

\begin{lemma}
\label{l:pf_thm_1}

 There exists a constant $c_{38}>0$ such that for all sufficiently large integer $n$, say $n\ge n_0$, we have
 $$
 \hbox{either } \prod_{i\in (\frac{n}{2}, \, n] \cap \z} [\E(m^{Y_i})]^{m-1} \ge 8,
 \quad
 \hbox{or}
 \quad
 \prod_{i=0}^n [\E(m^{Y_i})]^{m-1} \ge c_{38} \, n^{\alpha-2}\, .
 $$

\end{lemma}

\medskip

\noindent {\it Proof.} Let $M_0\ge 1$ and $c_{33}\ge 1$ be the constants in \eqref{cor1_variant_stable} in Remark \ref{r:Z2upperQ}. Let $c_{39}:= [120(m-1)c_{33}]^{2/(\alpha-2)} \ge 120$. Let $n\ge \lceil 3 c_{39} \, M_0 \rceil =: n_0$ be an integer, and let $M = M(n):= \lfloor \frac{n}{c_{39}} \rfloor$. So $\frac{n}{2c_{39}} \le M \le \frac{n}{120}$. Let $u_n := m- \frac{c_{40}}{n}$, where $c_{40}:= 30m$. We can enlarge the value of $M_0$ if necessary to ensure that $u_n > \frac{m}{2}$.

We discuss on two possible situations, each leading to one of the inequalities stated in the lemma.

{\it First situation: $\E[(1-(m-1)Y_i)u_n^{Y_i}] < \frac12$ for all $i\in (\frac{n}{2}, \, n] \cap \z$.} In this situation, we have, for all $i\in (\frac{n}{2}, \, n] \cap \z$, $(m-1) \E (Y_i \, u_n^{Y_i}) \ge \E (u_n^{Y_i}) - \frac12 \ge \frac12$, thus
$$
\E (Y_i \, u_n^{Y_i}) \ge \frac{1}{2(m-1)} \, .
$$

\noindent Consider the function $s\mapsto f_i(s) := \E(s^{Y_i})$, $s\in [0, \, m]$. We have
$$
\E(m^{Y_i})
=
f_i(m)
\ge
f_i(u_n) + (m-u_n)f_i'(u_n)
\ge
1+(m-u_n)f_i'(u_n) .
$$

\noindent Since $m-u_n= \frac{30m}{n}$ and $f_i'(u_n) = \frac{1}{u_n}\, \E (Y_i \, u_n^{Y_i}) \ge \frac1m \, \frac{1}{2(m-1)} = \frac{1}{2m(m-1)}$, we get, for all $i\in (\frac{n}{2}, \, n] \cap \z$,
$$
\E(m^{Y_i}) \ge 1+ \frac{15}{(m-1)n} \, .
$$

\noindent Consequently,
$$
\prod_{i\in (\frac{n}{2}, \, n] \cap \z} [\E(m^{Y_i})]^{m-1}
\ge
\Big( 1+ \frac{15}{(m-1)n} \Big)^{\! (m-1)n/2}
\ge 8 \, .
$$

{\it Second (and last) situation: $\E[(1-(m-1)Y_\ell)u_n^{Y_\ell}] \ge \frac12$ for some $\ell = \ell(n) \in (\frac{n}{2}, \, n] \cap \z$.} We will be working with this particular $\ell$ in the rest of the proof. Let $Z_0 := Y_0 \, {\bf 1}_{\{ Y_0 \le M\} }$ as in \eqref{Z0}. Let $(Z_n, \, n\ge 0)$ be a Derrida--Retaux system whose initial distribution is given by $Z_0$.

Since $u_n \in [1, \, m]$, the function $x\mapsto (1-(m-1)x)u_n^x$ is decreasing on $[0, \, \infty)$, so
\begin{equation}
    \E[(1-(m-1)Z_\ell)u_n^{Z_\ell}]
    \ge
    \E[(1-(m-1)Y_\ell)u_n^{Y_\ell}]
    \ge
    \frac12\, .
    \label{phi(un)}
\end{equation}

\noindent Consider the function
$$
\varphi(s) := \E[(1-(m-1)Z_\ell)s^{Z_\ell}] ,
\qquad
s\ge 0 ,
$$

\noindent which is well-defined because $Z_\ell$ is bounded. We have just proved that $\varphi(u_n) \ge \frac12$. Let $v_n := m+ \frac{1}{2c_4 n}$, where $c_4 \ge 1$ is the constant in Theorem $\ref{t:H_nkupper}$. Since $v_n > m$, we have, by concavity of $\varphi(\cdot)$,
$$
\varphi(v_n) \ge \varphi(m) + (v_n-m) \varphi'(v_n) \, .
$$

\noindent By assumption, the system $(Y_n, \, n\ge 0)$ is critical, so $(Z_n, \, n\ge 0)$ is subcritical (or critical in case $Y_0 \le M$ a.s.), which implies that $(m-1) \E(Z_\ell \, m^{Z_\ell}) \le \E(m^{Z_\ell})$. This means that $\varphi(m) \ge 0$. On the other hand, $v_n-m = \frac{1}{2c_4 n}$, whereas
$$
\varphi'(v_n)
=
\E[(1-(m-1)Z_\ell)Z_\ell v_n^{Z_\ell-1}]
\ge
-(m-1) \E(Z_\ell^2 v_n^{Z_\ell-1}),
$$

\noindent which is $= - \frac{m-1}{v_n} \E(Z_\ell^2 v_n^{Z_\ell}) \ge - \frac{m-1}{m} \E(Z_\ell^2 v_n^{Z_\ell})$. Assembling these pieces together yields that
$$
\varphi(v_n) \ge - \frac{1}{2c_4 n} \, \frac{m-1}{m} \, \E(Z_\ell^2 v_n^{Z_\ell})\, .
$$

\noindent Since $v_n = m+ \frac{1}{2c_4 n} \le m+ \frac{1}{2c_4 \ell} = m+ \frac{1}{2c_4 (\ell\vee M)}$ (to obtain the last equality, we have used the fact that $\ell\ge \frac{n}{2} \ge M$), we are entitled to apply inequality \eqref{cor1_variant_stable} in Remark \ref{r:Z2upperQ} to see that
$$
\E(Z_\ell^2 v_n^{Z_\ell})
\le
c_{33} (\ell\vee M)^{(4-\alpha)/2} M^{(\alpha-2)/2} \, .
$$

\noindent Since $(\ell \vee M) = \ell \le n$ and $M\le \frac{n}{c_{39}}$, this yields
$$
\E(Z_\ell^2 v_n^{Z_\ell})
\le
c_{33} n^{(4-\alpha)/2} (\frac{n}{c_{39}})^{(\alpha-2)/2}
=
\frac{c_{33}}{c_{39}^{(\alpha-2)/2}} \, n\, .
$$

\noindent Accordingly,
$$
\varphi(v_n)
\ge
- \frac{1}{2c_4 n} \, \frac{m-1}{m} \, \frac{c_{33}}{c_{39}^{(\alpha-2)/2}} \, n
=
- \frac{c_{41}}{c_{39}^{(\alpha-2)/2}} \, ,
$$

\noindent where $c_{41} :=\frac{c_{33}}{2c_4} \, \frac{m-1}{m}$. On the other hand, we have $\varphi(u_n) \ge \frac12$ (see \eqref{phi(un)}). Since $m= \beta u_n + (1-\beta) v_n$ with $\beta := \frac{1}{1+2c_4c_{40}} \in (0, \, 1)$, it follows from concavity of $\varphi(\cdot)$ that
$$
\varphi(m)
\ge
\beta \, \varphi(u_n) + (1-\beta) \, \varphi(v_n)
\ge
\frac{\beta}{2} - (1-\beta) \frac{c_{41}}{c_{39}^{(\alpha-2)/2}}
=
\frac{1- \frac{4c_4c_{40}c_{41}}{c_{39}^{(\alpha-2)/2}}}{2(1+2c_4c_{40})} \, .
$$

\noindent Our choice of the constant $c_{39}$ ensures $\frac{4c_4c_{40}c_{41}}{c_{39}^{(\alpha-2)/2}} = \frac12$; hence $\varphi(m) \ge \frac{1}{4(1+2c_4c_{40})} =: c_{42}$, i.e.,
$$
\E(m^{Z_\ell})- (m-1)\E(Z_\ell\, m^{Z_\ell})
\ge
c_{42} \, .
$$

\noindent Recall from \eqref{recursion_Coll} that
$$
\prod_{i=0}^{\ell-1} [\E(m^{Z_i})]^{m-1}
=
\frac{\E(m^{Z_\ell})- (m-1)\E(Z_\ell\, m^{Z_\ell})}{\E(m^{Z_0})- (m-1) \E(Z_0\, m^{Z_0})} \, .
$$

\noindent On the right-hand side, the numerator is at least $c_{42}$, whereas the denominator has already appeared in \eqref{E(mZ0)}:
\begin{eqnarray*}
    \E(m^{Z_0})- (m-1) \E(Z_0\, m^{Z_0})
 &=& \E[((m-1)Y_0 -1)\, m^{Y_0}\, {\bf 1}_{\{ Y_0 > M\} }) + \P(Y_0 > M)
    \\
 &\le& (m-1) \E(Y_0 \, m^{Y_0}\, {\bf 1}_{\{ Y_0 > M\} }) \, .
\end{eqnarray*}

\noindent Our assumption $\P(Y_0=j) \sim c_0 \, m^{-j}j^{-\alpha}$ (for $j\to \infty$) yields that $\E(Y_0 \, m^{Y_0}\, {\bf 1}_{\{ Y_0 > M\} }) \le \frac{c_{43}}{M^{\alpha-2}}$ for some constant $c_{43}>0$ depending on the law of $Y_0$. As a consequence,
$$
\prod_{i=0}^{\ell-1} [\E(m^{Z_i})]^{m-1}
\ge
\frac{c_{42}}{c_{43}\, M^{-(\alpha-2)}}
=
c_{44} \, M^{\alpha-2} \, ,
$$

\noindent where $c_{44}:= \frac{c_{42}}{c_{43}}$. This implies that
$$
\prod_{i=0}^{n-1} [\E(m^{Y_i})]^{m-1}
\ge
\prod_{i=0}^{\ell-1} [\E(m^{Y_i})]^{m-1}
\ge
\prod_{i=0}^{\ell-1} [\E(m^{Z_i})]^{m-1}
\ge
c_{44}\, M^{\alpha-2}\, .
$$

\noindent Since $M\ge \frac{n}{2c_{39}}$, this completes the proof of the lemma.\qed

\bigskip

We have all the ingredients for the proof of the lower bound in Theorem \ref{t:main}.

\bigskip

\noindent {\it Proof of the lower bound in Theorem \ref{t:main}.} It follows quite easily from Lemma \ref{l:pf_thm_1}. Indeed, let $n_0$ be the integer in Lemma \ref{l:pf_thm_1}, and let $n\ge 2n_0$. According to Lemma \ref{l:pf_thm_1}, there are two possibilities:

$\bullet$ either we have $\prod_{i\in (n/2^{j+1}, \, n/2^j] \cap \z} [\E(m^{Y_i})]^{m-1} \ge 8$ for all $0\le j \le \lfloor \frac{\log (n/n_0)}{\log 2}\rfloor$, in which case we have
$$
\prod_{i=0}^n [\E(m^{Y_i})]^{m-1}
\ge
8^{\lfloor \frac{\log (n/n_0)}{\log 2}\rfloor}
\ge
8^{\frac{\log (n/n_0)}{\log 2}-1}
=
\frac{n^3}{8n_0^3}\, ,
$$

\noindent which yields the lower bound in Theorem \ref{t:main} because $3>\alpha-2$.

$\bullet$ or there exists an integer $0\le j^* = j^*(n) \le \lfloor \frac{\log (n/n_0)}{\log 2}\rfloor$ such that $\prod_{i=0}^{n/2^{j^*}} [\E(m^{Y_i})]^{m-1} \ge c_{38} \, (\frac{n}{2^{j^*}})^{\alpha-2}$ and that $\prod_{i\in (n/2^{j+1}, \, n/2^j] \cap \z} [\E(m^{Y_i})]^{m-1} \ge 8$ for all non-negative integer $j<j^*$ (if there is any); in this case, we have (recalling notation: $\prod_\varnothing := 1$)
\begin{eqnarray*}
    \prod_{i=0}^n [\E(m^{Y_i})]^{m-1}
 &\ge& \prod_{i=0}^{n/2^{j^*}} [\E(m^{Y_i})]^{m-1} \times \prod_{j=0}^{j^*-1} \prod_{i\in (n/2^{j+1}, \, n/2^j] \cap \z} [\E(m^{Y_i})]^{m-1}
    \\
 &\ge& c_{38} \, (\frac{n}{2^{j^*}})^{\alpha-2} \times 8^{j^*} ,
\end{eqnarray*}

\noindent which implies that
$$
\prod_{i=0}^n [\E(m^{Y_i})]^{m-1}
\ge
c_{38} \, n^{\alpha-2}\, 2^{(5-\alpha)j^*}
\ge
c_{38} \, n^{\alpha-2}\, .
$$

In both situations, the lower bound in Theorem \ref{t:main} is valid.\qed

\bigskip
\bigskip

\noindent {\bf Acknowledgements.} We are grateful to Victor Dagard, Bernard Derrida, Yueyun Hu and Mikhail Lifshits for stimulating discussions on various aspects of the Derrida--Retaux model throughout the last few years. We also wish to thank two anonymous referees whose comments have led to improvements in the presentation of the paper.

\end{document}